\documentclass[a4paper,11pt,twoside,reqno]{amsart}

\usepackage[utf8]{inputenc}
\usepackage[plainpages=false,pdfpagelabels=true]{hyperref}
\usepackage{amssymb,amsthm}
\usepackage[margin=1in]{geometry}

\newtheorem{Satz}{Theorem}[section]

\newtheorem{Lem}[Satz]{Lemma}

\newtheorem{Cor}[Satz]{Corollary}

\newcommand{\tr}{{\operatorname{tr}}}
\theoremstyle{definition}

\newcommand{\dv}{\text{}dV}

\allowdisplaybreaks[1]

\renewcommand{\epsilon}{\varepsilon}

\newcommand{\R}{\ensuremath{\mathbb{R}}}

\numberwithin{equation}{section}

\title[A Liouville theorem for biharmonic maps between complete Riemannian manifolds]{A Liouville-type theorem for biharmonic maps between complete Riemannian manifolds with small energies}
\author{Volker Branding}
\date{\today}
\address{University of Vienna, Faculty of Mathematics\\
Oskar-Morgenstern-Platz 1, 1090 Vienna, Austria\\}
\email{volker.branding@univie.ac.at}
\subjclass[2010]{58E20; 53C43}
\keywords{biharmonic maps; complete Riemannian manifolds; Liouville theorem}
\begin{document}

\begin{abstract}
We prove a Liouville-type theorem for biharmonic maps from a complete Riemannian manifold of dimension \(n\) that has a lower bound
on its Ricci curvature and positive injectivity radius into a Riemannian manifold whose sectional curvature is bounded from above.
Under these geometric assumptions we show that if the $L^p$-norm of the tension field is bounded 
and the $n$-energy of the map is sufficiently small then every biharmonic map must be harmonic, where \(2<p<n\).
\end{abstract} 

\maketitle

\section{Introduction and Results}
\emph{Harmonic maps} belong to the most studied variational problems in differential geometry.
They are defined as critical points of the Dirichlet energy
\begin{align}
\label{harmonic-energy}
E_1(\phi)=\int_M|d\phi|^2\dv,
\end{align}
where \(\phi\colon M\to N\) is a map between the two Riemannian manifolds \((M,h)\) and \((N,g)\).
The critical points of \eqref{harmonic-energy} are characterized by the vanishing of the so-called
\emph{tension field}, which is given by
\begin{align*}
0=\tau(\phi):=\tr_h\nabla d\phi.
\end{align*}
This is a semilinear, elliptic second order partial differential equation, for which many results 
on existence and qualitative behavior of its solutions have been obtained.
For a recent survey on harmonic maps see \cite{MR2389639}.

A natural generalization of the harmonic map equation are the so-called \emph{biharmonic maps}.
These arise as critical points of the bienergy \cite{MR886529}, which is defined as
\begin{align*}
E_2(\phi)=\int_M|\tau(\phi)|^2\dv.
\end{align*}
In contrast to the harmonic map equation, the biharmonic map equation 
is of fourth order and is characterized by the vanishing of the \emph{bitension field}
\begin{align}
\label{euler-lagrange}
0=\tau_2(\phi):=-\Delta\tau(\phi)+R^N(d\phi(e_\alpha),\tau(\phi))d\phi(e_\alpha),
\end{align}
where \(e_\alpha\) is an orthonormal basis of \(TM\) and \(R^N\) denotes the curvature tensor
of the target manifold \(N\). We make use of the Einstein summation convention, meaning
that we sum over repeated indices.
The connection Laplacian on \(\phi^\ast TN\) is defined by
\(\Delta:=\nabla_{e_\alpha}\nabla_{e_\alpha}-\nabla_{\nabla_{e_\alpha}e_\alpha}\).
We choose the following convention for the 
Riemann curvature tensor \(R(X,Y)Z=[\nabla_X,\nabla_Y]Z-\nabla_{[X,Y]}Z\).

Note that the bienergy is conformally invariant in dimension \(\dim M=4\),
whereas the energy for harmonic maps is conformally invariant in dimension \(\dim M=2\).

As for harmonic maps there exists a large number of results on various aspects of biharmonic maps,
for a survey we refer to \cite{MR2301373}.

It is obvious that both harmonic maps and trivial maps solve the biharmonic map equation \eqref{euler-lagrange}.
For this reason it is interesting to find conditions under which solutions of the biharmonic map equation
reduce to harmonic or trivial maps. Several results in this direction have already been established:

\begin{enumerate}
 \item For biharmonic functions on flat space a Liouville theorem was already established before the first article on biharmonic maps appeared.
  Let \(M=\R^m\) and \(N=\R\) and consider a solution \(u\colon\R^m\to\R\) of \(\Delta^2u=0\).
  If \(u\) is bounded from above and below, then u is constant \cite{MR0274792}.
 \item Already in the first article on biharmonic maps \cite{MR886529} it was shown by a direct application of the maximum principle
  that biharmonic maps between compact Riemannian manifolds and the target having negative curvature must be harmonic.
 \item If the domain manifold is a complete non-compact Riemannian manifold with positive Ricci curvature
  and the target a Riemannian manifold of non-positive sectional curvature
  then every biharmonic map of finite bi-energy must be harmonic \cite{MR2604617}.
\item This result remains true if the condition of positive Ricci curvature on the domain is replaced by either the map having
finite energy or the domain manifold having finite volume \cite{MR3175248}.
\item This result was further generalized in \cite{MR3205803}: 
It is shown that a biharmonic map from a complete Riemannian manifold to a Riemannian manifold with non-positive
sectional curvature is harmonic if either the tension field is integrable in \(L^p\) with \(p\geq 2\)
and the volume of the domain is infinite or if the tension field is integrable in \(L^p\) with \(p\geq 2\)
and the Dirichlet energy of the map is finite.
\item Further generalizations were given in \cite{MR3427133} and \cite{MR3528391}.
\end{enumerate}

In this article we aim to establish another Liouville type theorem for biharmonic maps
that does not require any assumption on the curvature of the target.
To this end let \((M,h)\) and \((N,g)\) be two complete Riemannian manifolds.

We will not have to assume that the target manifold has negative curvature,
instead we will demand that the differential of the map is small in a certain \(L^p\)-norm.

In the following the positive constant \(A\) will
denote an upper bound on the sectional curvature \(K^N\) of \(N\).
We will denote the dimension of \(M\) by \(\dim M=n\).
The main result of this note is the following

\begin{Satz}
\label{theorem-liouville}
Let \(\phi\colon M\to N\) be a smooth biharmonic map.
Suppose that the Ricci curvature of \(M\) is bounded from below and that the injectivity radius of \(M\) is positive.
Assume that the sectional curvature \(K^N\) of \(N\) satisfies \(K^N\leq A\).
If 
\begin{align*}
\int_M|\tau(\phi)|^p\dv<\infty
\end{align*}
and 
\begin{align*}
\int_M|d\phi|^n\dv<\epsilon
\end{align*}
for \(2<p<n\) and \(\epsilon>0\) sufficiently small, then \(\phi\) must be harmonic.
\end{Satz}

If \(\phi\colon M\to N\) is an isometric immersion, we have \(\tau(\phi)=nH\),
where \(H\) denotes the mean curvature vector of the immersion.
The last theorem implies the following geometric result:
\begin{Cor}
Let \(\phi\colon M\to N\) be a smooth biharmonic isometric immersion.
Suppose that the Ricci curvature of \(M\) is bounded from below and that the injectivity radius of \(M\) is positive.
Assume that the sectional curvature \(K^N\) of \(N\) satisfies \(K^N\leq A\).

If
\begin{align*}
\int_M|H|^p\dv<\infty
\end{align*}
and 
\begin{align*}
\int_M|d\phi|^n\dv<\epsilon
\end{align*}
for \(2<p<n\) and \(\epsilon>0\) sufficiently small, then \(\phi\) must be a minimal immersion.
\end{Cor}

\section{Proof of the Theorem}
In order to prove Theorem \ref{theorem-liouville} we extend the ideas developed in \cite{MR3205803} and \cite{MR3175248}.
First, we need to establish some technical lemmas.

We choose a cutoff function  \(0\leq\eta\leq 1\) on \(M\) that satisfies
\begin{align*}
\eta(x)=1\textrm{ for } x\in B_r(x_0),\qquad \eta(x)=0\textrm{ for } x\in B_{2r}(x_0),\qquad |\nabla\eta|\leq\frac{C}{r}\textrm{ for } x\in M,
\end{align*}
where \(B_r(x_0)\) denotes the geodesic ball around \(x_0\) with radius \(r\).

\begin{Lem}
Let \(\phi\colon M\to N\) be a smooth biharmonic map.
Suppose that the sectional curvature \(K^N\) of \(N\) satisfies \(K^N\leq A\).
Then the following inequality holds
\begin{align}
\label{inequality-lp}
\frac{1}{2}\int_M\eta^2|\nabla\tau(\phi)|^2|\tau(\phi)|^{p-2}\dv\leq &A\int_M\eta^2|\tau(\phi)|^{p}|d\phi|^2\dv
+\frac{C}{R^2}\int_M|\tau(\phi)|^{p}\dv\\
\nonumber&-(p-2)\int_M\eta^2|\langle\nabla\tau(\phi),\tau(\phi)\rangle|^2|\tau(\phi)|^{p-4}\dv.
\end{align}
\end{Lem}

\begin{proof}
We test the biharmonic map equation \eqref{euler-lagrange} with \(\eta^2\tau(\phi)|\tau(\phi)|^{p-2}\) and find
\begin{align*}
\eta^2|\tau(\phi)|^{p-2}\langle\Delta\tau(\phi),\tau(\phi)\rangle=\eta^2|\tau(\phi)|^{p-2}\langle R^N(d\phi(e_\alpha),\tau(\phi))d\phi(e_\alpha),\tau(\phi)\rangle.
\end{align*}
Integrating over \(M\) and using integration by parts we obtain
\begin{align*}
\int_M\eta^2|\tau(\phi)|^{p-2}\langle\Delta\tau(\phi),\tau(\phi)\rangle\dv=&-2\int_M\langle\nabla\tau(\phi),\tau(\phi)\rangle|\tau(\phi)|^{p-2}\eta\nabla\eta\dv \\
&-(p-2)\int_M\eta^2|\langle\nabla\tau(\phi),\tau(\phi)\rangle|^2|\tau(\phi)|^{p-4}\dv \\
&-\int_M\eta^2|\nabla\tau(\phi)|^2|\tau(\phi)|^{p-2}\dv
\\
\leq&\frac{C}{R^2}\int_M|\tau(\phi)|^{p}\dv
-\frac{1}{2}\int_M\eta^2|\nabla\tau(\phi)|^2|\tau(\phi)|^{p-2}\dv\\
&-(p-2)\int_M\eta^2|\langle\nabla\tau(\phi),\tau(\phi)\rangle|^2|\tau(\phi)|^{p-4}\dv,
\end{align*}
where we used Young's inequality and the properties of the cutoff function \(\eta\).
Combining both equations we find
\begin{align*}
\frac{1}{2}\int_M\eta^2|\nabla\tau(\phi)|^2|\tau(\phi)|^{p-2}\dv\leq &\frac{C}{R^2}\int_M|\tau(\phi)|^{p}\dv
-(p-2)\int_M\eta^2|\langle\nabla\tau(\phi),\tau(\phi)\rangle|^2|\tau(\phi)|^{p}\dv \\
&-\int_M\eta^2|\tau(\phi)|^{p-2}\langle R^N(d\phi(e_\alpha),\tau(\phi))d\phi(e_\alpha),\tau(\phi)\rangle\dv.
\end{align*}
The result follows by estimating the last term on the right hand side.
\end{proof}

Now we recall the following fact:
Let \((M,g)\) be a complete Riemannian manifold whose Ricci curvature is bounded from below and 
with positive injectivity radius. Then for \(f\in W^{1,s}(M)\) with compact support the following Gagliardo-Nierenberg type inequality holds
\begin{align}
\label{gag-nie}
\|f\|_{L^r}\leq C\|df\|_{L^s},\qquad \frac{1}{r}=\frac{1}{s}-\frac{1}{n}, \qquad 1\leq s<n,
\end{align}
see \cite[Corollary 3.19]{MR1481970}.
This inequality allows us to give the following
  
\begin{Lem}
Let \((M,g)\) be a complete Riemannian manifold whose Ricci curvature is bounded from below and 
with positive injectivity radius. 
For \(p<n\) the following inequality holds
\begin{align}
\label{inequality-sobolev}
\int_M\eta^2|\tau(\phi)|^{p}|d\phi|^2\dv\leq C\big(\int_M|d\phi|^n\dv\big)^\frac{2}{n}\big(\frac{1}{R^2}\int_M|\tau(\phi)|^p\dv
+\int_M\eta^2|d|\tau(\phi)|^\frac{p}{2}|^2\dv\big),
\end{align}
where the positive constant \(C\) depends on \(n,p\) and the geometry of \(M\).
\end{Lem}

\begin{proof}
We set \(f:=|\tau(\phi)|^\frac{p}{2}\), which allows us to write
\begin{align*}
\int_M\eta^2|\tau(\phi)|^{p}|d\phi|^2\dv&=\int_M\eta^2f^{2}|d\phi|^2\dv.
\end{align*}
By Hölder's inequality we find
\begin{align*}
\int_M\eta^2f^{2}|d\phi|^2\dv\leq\big(\int_M(\eta f)^{2r}\dv\big)^\frac{1}{r}\big(\int_M|d\phi|^{\frac{2r}{r-1}}\dv\big)^{\frac{r-1}{r}}.
\end{align*}
Note that we may now apply the Gagliardo-Nierenberg type inequality \eqref{gag-nie} since \(\eta f\) has compact support.
Consequently, we find
\begin{align*}
\big(\int_M|\eta f|^{2r}\dv\big)^\frac{1}{2r}\leq C\big(\int_M|d(\eta f)|^\frac{2rn}{2r+n}\dv\big)^\frac{2r+n}{2rn}.
\end{align*}
Now, we choose \(r=\frac{n}{n-2}\) and obtain
\begin{align*}
\big(\int_M(\eta f)^{\frac{2n}{n-2}}\dv\big)^\frac{n-2}{n}\leq C\int_M|d(\eta f)|^2\dv.
\end{align*}
Using the properties of the cutoff function \(\eta\) we find
\begin{align*}
\int_M\eta^2f^{2}|d\phi|^2\dv\leq C\big(\int_M|d\phi|^n\dv\big)^\frac{2}{n}
\big(\int_M|d\eta|^2f^2\dv+\int_M\eta^2|df|^2\dv\big).
\end{align*}
The result follows by using that \(f=|\tau(\phi)|^\frac{p}{2}\).
\end{proof}

\begin{Lem}
Let \(\phi\colon M\to N\) be a smooth biharmonic map.
Suppose that the Ricci curvature of \(M\) is bounded from below and that the injectivity radius of \(M\) is positive.
Moreover, assume that the sectional curvature \(K^N\) of \(N\) satisfies \(K^N\leq A\).

If 
\begin{align*}
\int_M|\tau(\phi)|^p\dv <\infty
\end{align*}
and 
\begin{align*}
\int_M|d\phi|^n\dv<\epsilon
\end{align*}
for \(2<p<n\) and \(\epsilon>0\) sufficiently small, then \(|\tau(\phi)|\) is constant.
\end{Lem}
\begin{proof}
First of all, we note that
\begin{align*}
|\langle\nabla\tau(\phi),\tau(\phi)\rangle|^2|\tau(\phi)|^{p-4}=
\frac{1}{4}\big|d|\tau(\phi)|^2\big|^2|\tau(\phi)|^{p-4}
=\big|d|\tau(\phi)|\big|^2|\tau(\phi)|^{p-2}
=\frac{4}{p^2}\big|d|\tau(\phi)|^\frac{p}{2}\big|^2.
\end{align*}

Combining \eqref{inequality-lp} with \eqref{inequality-sobolev} and making use of the smallness
assumption of the energy of the map we find
\begin{align*}
\frac{1}{2}\int_M\eta^2|\nabla\tau(\phi)|^2|\tau(\phi)|^{p-2}\dv\leq&\frac{C}{R^2}\int_M|\tau(\phi)|^p\dv
+(C\epsilon-\frac{4(p-2)}{p^2})\int_M\eta^2|d|\tau(\phi)|^\frac{p}{2}|^2\dv.
\end{align*}
Choosing \(\epsilon\) sufficiently small, making use of the finiteness assumption
of the tension field and taking the limit \(R\to\infty\) we get
\begin{align*}
\int_M|\nabla\tau(\phi)|^2|\tau(\phi)|^{p-2}\dv=0,
\end{align*}
which implies that \(\tau(\phi)\) is a parallel vector field and thus has constant norm.
\end{proof}

In the following we will make use of the following result due to Gaffney \cite{MR0062490}:
\begin{Satz}
\label{gaffney}
Let \((M,h)\) be a complete Riemannian manifold. If a \(C^1\) one-form \(\omega\)
satisfies 
\begin{align*}
\int_M|\omega|\dv<\infty \qquad\text{ and }\qquad \int_M|\delta\omega|\dv<\infty
\end{align*}
or, equivalently, a \(C^1\) vector field \(X\) defined by \(\omega(Y)=h(X,Y)\),
satisfies
\begin{align*}
\int_M|X|\dv<\infty \qquad\text{ and }\qquad \int_M\operatorname{div}(X)\dv<\infty,
\end{align*}
then
\begin{align*}
\int_M(\delta\omega)\dv=\int_M\operatorname{div}(X)\dv=0.
\end{align*}
\end{Satz}

Making use of Gaffney's result we can now prove Theorem \ref{theorem-liouville}.

\begin{proof}[Proof of Theorem \ref{theorem-liouville}]
We define a one-form \(\omega\) as follows
\begin{align*}
\omega(X):=|\tau(\phi)|^{\frac{p(n-1)}{n}-1}\langle d\phi(X),\tau(\phi)\rangle.
\end{align*}
Note that 
\begin{align*}
\int_M|\omega|\dv\leq\int_M|\tau(\phi)|^{\frac{p(n-1)}{n}}|d\phi|\dv\leq\big(\int_M|d\phi|^n\dv\big)^\frac{1}{n}\big(\int|\tau(\phi)|^p\dv\big)^{\frac{n-1}{n}}<\infty,
\end{align*}
where we made use of the finiteness assumptions of both energy and bienergy of the map \(\phi\).
We fix an orthonormal basis of \(TM\) denoted by \(e_\alpha,\alpha=1,\ldots,\dim M\) and compute
\begin{align*}
-\delta\omega=&\nabla_{e_\alpha}\omega(e_\alpha)\\
=&\nabla_{e_\alpha}(|\tau(\phi)|^{\frac{p(n-1)}{n}-1}\langle d\phi(e_\alpha),\tau(\phi)\rangle)\\
=&|\tau(\phi)|^{\frac{p(n-1)}{n}-1}\langle\nabla_{e_\alpha}d\phi(e_\alpha),\tau(\phi)\rangle\\
=&|\tau(\phi)|^{\frac{p(n-1)}{n}+1},
\end{align*}
where we used that \(\nabla_X\tau(\phi)=0\) for all vector fields \(X\) and that \(|\tau(\phi)|\) is constant.
Due to the assumptions we have that \(\frac{p(n-1)}{n}+1>0\).
Again, since \(|\tau(\phi)|\) is constant and \(\int_M|\tau(\phi)|^p\dv\) is bounded by assumption
we can apply Theorem \ref{gaffney} and conclude that \(\tau(\phi)=0\).
\end{proof}

\par\medskip
\textbf{Acknowledgements:}
The author gratefully acknowledges the support of the Austrian Science Fund (FWF) 
through the project P30749-N35 ``Geometric variational problems from string theory''.
\bibliographystyle{plain}
\bibliography{mybib}
\end{document}